\documentclass[journal]{IEEEtran}  % Comment this line out
\usepackage{tikz} 
\usepackage{amsmath} 
\usepackage{amssymb} 
\usepackage{algorithm}
\usepackage{enumerate}
\usepackage[noend]{algpseudocode}
\usepackage{graphicx}
\usepackage[utf8]{inputenc}
\usepackage[noadjust]{cite}
\usepackage[caption=false, font=footnotesize]{subfig}
\usepackage{todonotes}
\usepackage{pgfplots}

\title{A Dual Active-Set Solver for Embedded Quadratic Programming Using Recursive LDL$^T$ Updates}
\author{Daniel Arnström, Alberto Bemporad and Daniel Axehill}
\begin{document}
\definecolor{set19c1}{HTML}{E41A1C}
\definecolor{set19c2}{HTML}{377EB8}
\definecolor{set19c3}{HTML}{4DAF4A}
\definecolor{set19c4}{HTML}{984EA3}
\definecolor{set19c5}{HTML}{FF7F00}
\definecolor{set19c6}{HTML}{FFFF33}
\definecolor{set19c7}{HTML}{A65628}
\definecolor{set19c8}{HTML}{F781BF}
\definecolor{set19c9}{HTML}{999999}
\maketitle
\newtheorem{proposition}{Proposition}
\newtheorem{lemma}{Lemma}
\newtheorem{corollary}{Corollary}
\newtheorem{theorem}{Theorem}
\newtheorem{remark}{Remark}
\begin{abstract}
  In this paper we present a dual active-set solver for quadratic programming which has properties suitable for use in embedded model predictive control applications. In particular, the solver is efficient, can easily be warm-started, and is simple to code. Moreover, the exact worst-case computational complexity of the solver can be determined offline and, by using outer proximal-point iterations, ill-conditioned problems can be handled in a robust manner.  

\end{abstract}
\section{Introduction}
Efficient, reliable, and simple quadratic programming (QP) solvers are essential when model predictive control (MPC) is used on embedded systems in real-time applications, where a QP has to be solved at each time step under real-time constraints with limited memory and computational resources. 

Popular methods for solving these QPs are active-set methods \cite{bemporad2015quadratic,ferreau2014qpoases,schmid1994quadratic,saraf2019bounded}, interior-point methods \cite{wang2010fast,frison2020hpipm}, gradient projection methods \cite{10.1109CDC.2008.4738961,patrinos2014accelerated,richter2009real,giselsson2014improved}, and operator splitting methods \cite{stellato2020osqp}.

In particular, the active-set method in \cite{bemporad2015quadratic} (QPNNLS), which is based on reformulating the QP as a nonnegative least-squares (NNLS) problem, is simple to implement and has proven to be efficient for solving small to medium size QP problems. Furthermore, its reliability has been improved greatly in \cite{bemporad2017numerically} where outer proximal-point iterations are used to improve its numerical stability, and in \cite{9104668}, where QPNNLS is shown to be closely related to a primal active-set QP method applied to the dual problem, allowing the complexity certification method in \cite{journal2020} to be used to determine the exact computational complexity of QPNNLS.

In this paper we use insights from \cite{9104668} to propose a dual active-set method for quadratic programming which retains the favorable properties of QPNNLS (efficiency and simplicity) by making recursive updates to an LDL$^T$ factorization. In addition to retaining favorable properties, we show that operating directly on the dual QP instead of the NNLS reformulation used in \cite{bemporad2015quadratic} yields additional improvements: 
(i) Direct reusability of matrix factors when the linear term in the objective function and the constant term in the constraints change, which is relevant for MPC and when the active-set method  is combined with outer proximal-point iterations. 
(ii) Improved numerical stability already in the setup without proximal-point iterations. 
(iii) Improved efficiency by reducing intermediate computations stemming from the NNLS reformulation. 

The main contribution of this paper is, hence, showing how the properties of QPNNLS can be retained and improved by directly operating on the dual QP problem instead of an NNLS problem. Concretely this requires, for example, extending the recursive LDL$^T$ updates to detect and handle singularities.
More broadly, we show that a simple dual active-set solver that performs such recursive LDL$^T$ updates can outperform state-of-the-art solvers on small to medium size QPs,  which are often encountered in embedded MPC applications. 

The rest of the paper is structured as follows:
in the remainder of this introduction we introduce the quadratic programming problem and notation that is used throughout the paper. In Section \ref{sec:qpalg} we present the proposed dual active-set algorithm. We then discuss the relationship between this algorithm and other, similar, QP algorithms in Section \ref{sec:relate}. Extensions to the algorithm, addressing some practical concerns, are presented in Section \ref{sec:extensions}. 
Finally, we present numerical experiments in Section \ref{sec:result}. 

\subsection{Problem formulation}
We consider a quadratic program in the form 
\begin{equation}
  \label{eq:QP}
  \begin{aligned}
	&\underset{x}{\text{minimize}}&&J(x) \triangleq\frac{1}{2}x^T H x + f^T x \\
	&\text{subject to} &&Ax \leq b,\\
  \end{aligned}
\end{equation}
where $x\in \mathbb{R}^n$. The objective function $J$ is given by $H\in \mathbb{S}^n_{++}$ and $f\in \mathbb{R}^n$, and the feasible set is given by $A\in \mathbb{R}^{m\times n}$ and $b\in \mathbb{R}^m$. Furthermore, the solution to \eqref{eq:QP} is denoted by $x^*$.  
\begin{remark}[Relaxing strict convexity]
  The proposed algorithm requires $H$ to be invertible. The case when $H\succeq0$ can, however, be handled by performing proximal-point iterations, described in Section \ref{ssec:proximal}.
\end{remark}

  A set of necessary and, because of the convexity of \eqref{eq:QP}, sufficient conditions for optimality of $x^*$ are the KKT-conditions:%
\begin{subequations}
  \label{eq:kkt}
\begin{align}
  H x^* + A^T \lambda^*  &= -f, \label{eq:kktstat} \\
  A x^* \leq b,\: \lambda^*&\geq 0, \\
  [b-A x^*]_i [\lambda^*]_i &= 0,\:\: \forall i=1,\dots, m,
\end{align}
for $\lambda^* \in \mathbb{R}^m$, and where the operator $[\cdot]_i$ extracts the $i$th row of a matrix, or the $i$th entry of a vector.
\end{subequations}

Instead of solving \eqref{eq:QP} directly, we will solve its so-called \textit{dual problem}: 
\begin{equation}
  \label{eq:qp-dual}
	 \underset{\lambda \geq 0}{\text{minimize }} J_d(\lambda) \triangleq \frac{1}{2}\lambda^T M M^T \lambda + d^T \lambda,
    % &\text{subject to} && \lambda \geq 0.
\end{equation}
where we have, similar to \cite{bemporad2015quadratic}, introduced 
\begin{equation}
  \label{eq:aux-def}
  M \triangleq A R^{-1}, \quad v\triangleq R^{-T} f,\quad d\triangleq b+M v,
\end{equation}
and where $R$ is an upper triangular Cholesky factor of $H$ ${(H = R^T R)}$. 
The solution $\lambda^*$ to \eqref{eq:qp-dual} satisfies the same KKT-conditions as \eqref{eq:QP} (see, e.g., \cite{dorn1960duality}) and $x^*$ can, hence, be recovered from $\lambda^*$ through \eqref{eq:kktstat} when $H \succ 0$. 
\subsection{Notation}
Since the proposed method, soon to be introduced, is an iterative method, $\lambda_k$ denotes the value of the dual iterate at iteration $k$. Furthermore, $\mathcal{W}_k$ is the so-called \textit{working set} which contains indices of the components of $\lambda_k$ which are free to vary (which can be interpreted as imposing the corresponding primal constraints to hold with equality). Conversely, $\overline{\mathcal{W}}_k$ contains the components of $\lambda_k$ which are fixed at zero. 
$M_k$ and $d_k$ denote the rows of $M$ and $d$ indexed by $\mathcal{W}_k$, respectively. Likewise, $\overline{M}_k$ and $\overline{d}_k$ denote the rows of $M$ and $d$ indexed by $\overline{\mathcal{W}}_k$, respectively. 
Finally, $\text{ker}(A)$ denotes the kernel of a matrix $A$.

\section{A Dual Active-Set Algorithm}
\label{sec:qpalg}
The dual active-set algorithm that we propose, given in Algorithm \ref{alg:dual-as}, can be interpreted as the primal active-set algorithm considered in \cite{journal2020} applied to the dual problem in \eqref{eq:qp-dual}, which is mathematically equivalent to several other popular active-set algorithm formulations (see Section \ref{ssec:other-as} for details and advantages of the proposed formulation).
We will now give an overview of the algorithm and then cover specifics, such as how to efficiently solve the subproblems encountered in Step \ref{step:nonsingular-subproblem}  and \ref{step:singular-subproblem}, in Section \ref{sec:ldl}.  
\subsection{Algorithm overview}
\label{ssec:overview}
Like any other active-set method, Algorithm \ref{alg:dual-as} iteratively updates the working set $\mathcal{W}$. An iteration always starts by solving an equality constrained subproblem defined by the current working set $\mathcal{W}_k$: 
\begin{equation}
  \label{eq:subproblem}
  \begin{aligned}
	 &\underset{\lambda}{\text{min}}&& \frac{1}{2}\lambda^T M M^T \lambda + d^T \lambda\\
	 &\text{s.t.} &&[\lambda]_i = 0,\: \forall i\notin \mathcal{W}_k,
  \end{aligned}
\end{equation}
where $k$ is the current iteration. By using the constraint $[\lambda]_{\overline{\mathcal{W}}_k}=0$ to eliminate variables, \eqref{eq:subproblem} is equivalent to the unconstrained problem 
\begin{equation}
  \label{eq:unconstrained}
  \underset{[\lambda]_{\mathcal{W}_k}}{\text{min}} \frac{1}{2}[\lambda]^T_{\mathcal{W}_k} M_k M_k^T [\lambda]_{\mathcal{W}_k} + d_k^T [\lambda]_{\mathcal{W}_k},\quad [\lambda]_{\overline{\mathcal{W}}_k}=0.
\end{equation}
If $M_k M_k^T \succ 0$ , the unconstrained problem in \eqref{eq:unconstrained} has a unique solution and the solution $\lambda^*_k$ to \eqref{eq:subproblem} is then given by 
\begin{equation}
  \label{eq:linsystem}
  M_k M_k^T [\lambda^*_k]_{\mathcal{W}_k} = -d_k, \quad [\lambda^*_k]_{\overline{\mathcal{W}}_k}=0.
\end{equation}
If $\lambda^*_k \geq 0$, we set $\lambda_{k+1} \leftarrow \lambda^*_k$ and check for primal feasibility (see below). Otherwise, a line-search along the line-segment connecting $\lambda_k$ and $\lambda^*_k$ is performed and the first component which becomes zero is removed from $\mathcal{W}_k$, i.e., is fixed at zero.

If $M_k M_k^T$ is singular and $d_k^T p \neq 0$ for some ${p\in \text{ker}(M_k M_k^T)}$, there is no bounded solution to \eqref{eq:unconstrained}, i.e., the objective function can be made arbitrarily small by moving in a direction $p_k$ which satisfies 
\begin{equation}
  \label{eq:singdireq}
  [p_k]_{\overline{\mathcal{W}}_k}=0,\quad M_k M_k^T [p_k]_{\mathcal{W}_k}=0,\quad d^T p_k < 0. 
\end{equation}
Hence, a line-search along the ray $\lambda_k + \alpha p_k, \alpha>0$ is performed. As is discussed in detail in Section \ref{ssec:infeas}, at least one component of $\lambda$ not in $\mathcal{W}_k$ will become zero while moving along this ray if \eqref{eq:QP} is feasible. Again, the first component which becomes zero is removed from $\mathcal{W}_k$, i.e., is fixed at zero.    

\begin{remark}[Impossibility of $d^T p_k = 0$]
  For Algorithm \ref{alg:dual-as}, one can show that once $M_k M_k^T$ becomes singular, there always exists a solution to \eqref{eq:singdireq}, see, e.g., Lemma 3.5 in \cite{arnstrom2021complexity} for details.
\end{remark}

When $\lambda^*_k \geq 0$, primal feasibility for the constraints not in $\mathcal{W}_k$ is checked by computing the primal slack $\mu_k$ (which is the dual vector of \eqref{eq:qp-dual})
\begin{equation}
  \label{eq:primal-slack}
[{\mu}_k]_{\overline{\mathcal{W}}_k} = \overline{M}_k M_k \lambda^*_k + \overline{d}_k.
\end{equation}
Primal feasibility is satisfied if ${\mu}_k \geq 0$ and, since stationarity, dual feasibility, and complimentary slackness already hold, $\lambda^*_k$ is optimal. Otherwise, the most negative component of $\mu_k$ is added to $\mathcal{W}_k$ (making the corresponding component of $\lambda$ free to vary).

\begin{remark}[Selection rule]
  Adding \textit{the most} negative component of $\mu_k$ to $\mathcal{W}_k$ is a common rule in practice, but adding \textit{any} negative component of $\mu_k$ to $\mathcal{W}_k$ also leads to convergence.
\end{remark}

\begin{remark}[Primal feasibility tolerance]
  When implemented in practice, $\mu_k \geq -\epsilon_p$ is considered instead of $\mu_k \geq 0$ in Step \ref{step:primalfeas} for numerical reasons, where $\epsilon_p >0$ is the tolerance for primal feasibility. Similar tolerances should also be used for the inequalities in Steps \ref{step:csp}, \ref{step:innerstart}, and \ref{step:singblocking}.
\end{remark}

After $\mathcal{W}_k$ has been changed by either adding or removing an index to get a new working set $\mathcal{W}_{k+1}$, the algorithm starts another iteration by solving \eqref{eq:unconstrained} for $\mathcal{W}_{k+1}$ (or by solving \eqref{eq:singdireq} if $M_{k+1} M_{k+1}^T$ is singular) and the steps described above are repeated until convergence. 
%\begin{remark}[Retreiving primal solution]
%  In the context of linear MPC, the entire primal solution $x^*$ does not have to be recovered from the dual solution $\lambda^*$ in Step \ref{step:return}, only the components corresponding to the control input at the current time step have to be computed. 
%\end{remark}

\begin{algorithm}[H]
  \caption{Dual active-set method for solving \eqref{eq:QP}.}
  \label{alg:dual-as}
  \begin{algorithmic}[1]
	\Require $M,d,v, R^{-1}, \mathcal{W}_0, \lambda_0$
	\Ensure $x^*, \lambda^*, \mathcal{A}^*$
	\algrenewcommand\algorithmicindent{0.825em}%
	\While{true}
	\If{$M_k M_k^T$ is nonsingular}\label{step:nonsingular-start}
	\State [$\lambda^*_k]_{\mathcal{W}_k} \leftarrow$ solution to $M_k M_k^T [\lambda^*_k]_{\mathcal{W}_k} = -d_k$ \label{step:nonsingular-subproblem}
	\If{$\lambda^*_k \geq 0$} \label{step:csp}
	\State $[\mu_k]_{\overline{\mathcal{W}}_k} \leftarrow \overline{M}_k M_k^T [\lambda^*_k]_{\mathcal{W}_k} + \overline{d}_k$,\: $\lambda_{k+1} \leftarrow \lambda^*_k$
	\If{$\mu_k \geq 0$} optimum found, \textbf{goto} \ref{step:return}\label{step:primalfeas}
	\Else \quad$j \leftarrow \text{argmin}_{i\in \overline{\mathcal{W}}_k}[\mu_k]_i$, \:\:$\mathcal{W}_{k+1} \leftarrow \mathcal{W}_k\cup\{j\}$
	\EndIf
	\Else
	%\Comment{$\exists i: [\lambda^*_k]_i < 0$ }
	\State $p_k \leftarrow \lambda^*_k-\lambda_k$, \:\:$\mathcal{B}\leftarrow\{i\in\mathcal{W}_k: [\lambda^*_k]_i < 0\}$ \label{step:innerstart} 
	\State $[\lambda_{k+1}, \mathcal{W}_{k+1}] \leftarrow$ \textsc{fixComponent}$(\lambda_k,\mathcal{W}_k,\mathcal{B},p_k)$ \label{step:innerstop}
	\EndIf
	\Else $\:(M_k M_k^T \text{singular})$\label{step:singular-start}
	%\Comment{$M_k M_k^T$ is singular}
	\State $[p_k]_{\mathcal{W}_k} \leftarrow$ solution to $M_k M^T_k [p_k]_{\mathcal{W}_k} =0,\:\: p_k^T d < 0$ \label{step:singular-subproblem}
	\State $\mathcal{B} \leftarrow \{i\in \mathcal{W}_k: [p_k]_i < 0\}$ \label{step:singblocking}
	\State $[\lambda_{k+1}, \mathcal{W}_{k+1}] \leftarrow$ \textsc{fixComponent}$(\lambda_k,\mathcal{W}_k,\mathcal{B},p_k)$ \label{step:singular-stop}
	\EndIf
	\State $k \leftarrow k+1$ 
	\EndWhile
	\State \Return $x^* \leftarrow  -R^{-1}(M^T_k [\lambda^*_k]_{\mathcal{W}_k}+v), \lambda^*_k, \mathcal{W}_k$ \label{step:return}
	\hrule\rule{0pt}{1pt}
	\Procedure{fixComponent}{$\lambda_k, \mathcal{W}_k,\mathcal{B},p_k$}
	\State $j \leftarrow \text{argmin}_{i \in \mathcal{B}} -[\lambda_k]_i/[p_k]_i $ \label{step:ratio} 
	\State $\mathcal{W}_{k+1} \leftarrow \mathcal{W}_k\setminus\{j\}$,\:\: $\lambda_{k+1} \leftarrow \lambda_k - ([\lambda_k]_j/[p_k]_j) p_k$
	\EndProcedure
  \end{algorithmic}
\end{algorithm}

The convergence of Algorithm \ref{alg:dual-as} can be proven by standard arguments for active-set methods (cf. e.g., Section 4 in \cite{fletcher1971general} or Section 3 in \cite{10.1007BF02591962}). 

\begin{remark}[Equality constraints]
  If equality constraints are included in \eqref{eq:QP}, Algorithm \ref{alg:dual-as} can be extended by, essentially, adding additional components to $\lambda$ for each added equality constraints and conceptually treating these as always being in $\mathcal{W}_k$, i.e, always being free (see Section \ref{ssec:eqc} for details). 
\end{remark}

\subsection{LDL$^T$ factorization}
\label{sec:ldl}
Matrix $M_k M_k^T$ is central in Algorithm \ref{alg:dual-as}, partly because of whether it is singular or not results in different modes (Steps \ref{step:nonsingular-start}-\ref{step:innerstop} or Steps \ref{step:singular-start}-\ref{step:singular-stop}, respectively), but also since it is used to compute $\lambda^*_k$ in Step \ref{step:nonsingular-subproblem} or $p_k$ in Step \ref{step:singular-subproblem}. We will now show that the above-mentioned operations can be efficiently performed by factorizing $M_k M_k^T=L D L^T$, where $L$ is a lower unit triangular matrix and $D$ is a diagonal matrix. In particular, we show in Section \ref{ssec:ldl-kkt-indentify} that the singularity of $M_k M_k^T$ can easily be indentified, and in Section \ref{ssec:ldl-kkt-nonsing} and Section \ref{ssec:ldl-kkt-sing} we show that the system of linear equations defining $\lambda^*_k$ and $p_k$ can be efficiently solved. Moreover, since only a single row of $M_k$ is either added or removed between iterations in Algorithm \ref{alg:dual-as}, $L$ and $D$ can be recursively updated,  as described in Section \ref{ssec:ldlupdate}, which reduces the computational complexity of the algorithm significantly.
\subsubsection{Detecting singularity}
\label{ssec:ldl-kkt-indentify}
Given an $L$ and $D$ it is straightforward to determine whether $M_k M_k^T$ is singular or not since this can directly be seen in $D$:
\begin{lemma}
 \label{lem:singular}
 Assume that there exist $L$ and $D$ such that $M_k M_k^T = LDL^T$, with $L$ being a lower unit triangular matrix and $D$ being a diagonal matrix with nonnegative elements. Then $[D]_{ii} \neq 0, \forall i \Leftrightarrow M_k M_k^T$ is nonsingular  
\end{lemma}

\begin{IEEEproof}
  Directly follows from $L$ having full rank and that $D$ is a diagonal matrix. 
\end{IEEEproof}
%\begin{IEEEproof}
%  Recall that for a triangular matrix $T\in \mathbb{R}^{N\times N}$,  $\det(T) = \prod_{i=1}^{N} [T]_{ii}$. Also recall the identity $\det(A B)=\det(A)\det(B)$. From this we get $\det(M_k M_k^T) = \det(L)\det(D)\det(L^T) = \prod_{i=1}^{N} [D]_{ii}$, where we have used that $\det(L)=\det(L^T)=1$ since $L$ is unit triangular. From this equality it is evident that $\det(M_k M_k)^T \neq 0 \Leftrightarrow [D]_{ii} \neq 0, \forall i$.
%\end{IEEEproof}
\subsubsection{Solving nonsingular KKT-systems}
\label{ssec:ldl-kkt-nonsing}
When $M_k M_k^T$ is nonsingular, solving ${M_k M_k^T [\lambda^*_k]_{\mathcal{W}_k} = -d_k}$ given an LDL$^T$ factorization of $M_k M_k^T$ can be done by solving two triangular linear systems: 
\begin{subequations}
  \label{eq:backsolve}
\begin{align}
  y &\leftarrow \text{Solve }L y = -d_k  \label{eq:kkt-solve-fwd} \text{ (forward substitution)}, \\ 
  z &\leftarrow \text{Scale } y \text{ with } D,\quad  [z]_i = \tfrac{[y]_i}{[D]_{ii}}\label{eq:kkt-solve-scale}, \\
  [\lambda^*_k]_{\mathcal{W}_k} &\leftarrow \text{Solve } L^T [\lambda^*_k]_{\mathcal{W}_k} = z \text{ (backward substitution)}.
\end{align}
\end{subequations}
\begin{remark}[Efficient forward substitution] Because a certain number of first rows of the lower-triangular matrix $L$ and of $d_k$ and $D$ remain constant between iterations, both the forward substitution and scaling in \eqref{eq:kkt-solve-fwd} and \eqref{eq:kkt-solve-scale} do not have to be done from scratch at each iteration of Algorithm \ref{alg:dual-as}. The amount of previous computations that can be reused depends on how much the working set $\mathcal{W}_k$ changes. For example, only the last element of $y$ in \eqref{eq:kkt-solve-fwd} has to be solved for after a constraint is added to $\mathcal{W}_k$.
\end{remark}
\subsubsection{Solving singular KKT-systems}
\label{ssec:ldl-kkt-sing}
Next, we consider the case when $M_k M_k^T$ is singular which, from Lemma \ref{lem:singular}, means that at least one diagonal element of $D$ is zero. The following lemma shows that the LDL$^T$ factorization can be used to efficiently compute a $p$ which satisfies $M_k M_k^T p = 0$. 
\begin{lemma} Assume that $M_k M_k^T = LDL^T$ with the $i$th diagonal element of $D$ being zero, i.e.,
\begin{equation}
  L = 
  \begin{bmatrix}
	L_1& 0& 0 \\
	l^T_i& 1 & 0 \\
	*& * & * 
  \end{bmatrix}, \quad
  D =
  \begin{bmatrix}
	D_1 & 0 & 0\\
	0 & 0 & 0 \\
	0 & 0 & D_2
  \end{bmatrix}.
\end{equation}
Let $\tilde{p}$ be the solution to $L_1^T \tilde{p} = -l_i$. Then $p =\begin{bmatrix}
  \tilde{p}^T & 1&  0 
\end{bmatrix}^T$ satifies $M_k M_k^T p=0$. 
\end{lemma}
\begin{IEEEproof}
  By multiplying $L^T$ with the given $p$ we get
\begin{equation}
  \label{eq:unit}
  L^T p = 
  \begin{bmatrix}
	L_1^T& l_i& *\\ 
	0& 1 & *\\ 
	0& 0 & * 
  \end{bmatrix}
  \begin{bmatrix}
	\tilde{p} \\ 1\\ 0 
  \end{bmatrix} 
  =
  \begin{bmatrix}
	L_1^T \tilde{p} + l_i \\ 1\\ 0 
  \end{bmatrix}
  = e_i
\end{equation}
where $e_i$ is the $i$th unit vector and we have used that $L_1^T \tilde{p} + l_i=0$. Using \eqref{eq:unit} and that the $i$th element of $D$ is zero gives 
\begin{equation}
  M_k M_k^T p = L D L^T p = L D e_i = 0, 
\end{equation}
proving the lemma.  
\end{IEEEproof}
Hence, we can find a nontrivial null vector of $M_k M_k^T$ by solving an $i$th dimensional upper unit triangular system of linear equations by backward substitution. In most cases, $i=|\mathcal{W}_k|$ since $M_k M_k^T$ only becomes singular after an addition of an element to $\mathcal{W}_k$, which in turn implies that the zero element in $D$ will be the last diagonal element (if the updates described below are used). 

How null vectors of $M_k M_k^T$ can be used to detect primal infeasibility is discussed in Section \ref{ssec:infeas}.
%If $M_k M_k^T$ remains singular after a constraint is subsequently removed the zero element might not, any longer, be the last diagonal element of $D$. 
%\todo[inline]{It will always be the last element which is singular, after a removal you cannot maintain singularity.}
\subsubsection{Updating LDL$^T$ after addition to/removal from $\mathcal{W}$}
\label{ssec:ldlupdate}
To recursively update $L$ and $D$ after adding an index to $\mathcal{W}$, we recall the result from Theorem 2 in \cite{bemporad2015quadratic}. 
\begin{lemma}
  Let $L$ be a unit lower triangular matrix and $D$ be a diagonal matrix such that $M_k M_k^T = L D L^T$. Furthermore let $M^{\text{+}} = [\begin{smallmatrix} 
	  M_k \\ [M]_i 
\end{smallmatrix}]$. Then $M^{\text{+}} (M^{\text{+}})^T = L^{\text{+}} D^{\text{+}} (L^{\text{+}})^T$ with 
  \begin{equation}
	L^{\text{+}} = 
  \begin{bmatrix}
	L &0\\
	l^T & 1
  \end{bmatrix}, \quad
  D^{\text{+}} =
  \begin{bmatrix}
	D &0 \\
	0 & \delta
  \end{bmatrix}, \\
  \end{equation}
  where $l$ and $\delta$ are defined by
  \begin{equation}
  L D l = M_k [M]_i^T, \quad \delta = [M]_i [M]_i^T - l^T D l.
  \end{equation}
\end{lemma}
\begin{IEEEproof}
  Cf. proof of Theorem 1 in \cite{bemporad2015quadratic}.
\end{IEEEproof}

Similarly, we recall the result from Lemma 2 in \cite{bemporad2015quadratic} for recursively updating $L$ and $D$ after the removal of an index from $\mathcal{W}$.

\begin{lemma}
  \label{lem:remove}
Let $M_k M_k^T = L D L^T$ with 
\begin{equation}
  L = 
  \begin{bmatrix}
	L_1& 0& 0\\
	*&  1& 0 \\
	A & l_i & L_2
  \end{bmatrix},\quad
  D = 
  \begin{bmatrix}
	D_1 & 0 &0 \\ 
	0& \delta_i & 0 \\
	0& 0& D_2
  \end{bmatrix},
\end{equation}
where $l_i$ and $\delta_i$ are placed in the $i$th column of $L$ and $D$, respectively, and where $L_1$ and $L_2$ are lower unit triangular and $D_1$ and $D_2$ are diagonal. Furthermore, let $M^{\text{-}}$ be $M_k$ with the $i$th row removed. Then $M^{\text{-}} (M^{\text{-}})^T = L^{\text{-}} D^{\text{-}} (L^{\text{-}})^T$, where $L^\text{-}$ and $D^{\text{-}}$ are given by 
\begin{equation}
  L^{\text{-}} = 
  \begin{bmatrix}
	L_1&  0\\
	A & \tilde{L}_2
  \end{bmatrix},\quad
  D^{\text{-}} = 
  \begin{bmatrix}
	D_1  &0 \\ 
	0&  \tilde{D}_2
  \end{bmatrix},
\end{equation}
if $\tilde{L}_2 \tilde{D}_2 \tilde{L}_2^T = L_2 D_2 L_2^T +\delta_i l_i l_i^T$, where $\tilde{L}_2$ is lower unit triangular and $\tilde{D}_2$ is diagonal. 
\end{lemma}
\begin{IEEEproof}
  Cf. the proof of Lemma 2 in \cite{bemporad2015quadratic}.
\end{IEEEproof}
In Lemma \ref{lem:remove}, the LDL$^T$ factorization of $L_2 D_2 L_2^T + \delta_i l_i l_i^T$ is a rank-one update of an existing LDL$^T$ factorization, which can be done efficiently by, e.g., Algorithm C1 in \cite{gill1974methods}. 
\section{Comparison to similar QP algorithms}
\label{sec:relate}
\subsection{QPNNLS}
\label{ssec:qpnnlscomp}
The proposed method is similar to the QP method described in \cite{bemporad2015quadratic} (QPNNLS) which is based on transforming the QP to a nonnegative least-squares (NNLS) problem. This transformation results in solving least-squares subproblems on the form ${\min_{y}\|A_k^T y-b_k\|_2^2}$ at each iteration,  which in turn results in solving the linear equation 
  $A_k A_k^T y = A_k b_k$, similar to the linear equation \eqref{eq:linsystem} solved in Algorithm \ref{alg:dual-as}. Similar to the proposed algorithm, an LDL$^T$ factorization of $A_k A_k^T$ is used to efficiently solve these linear equations. Explicitly, $A_k = [M_k  \:\: d_k]$ in \cite{bemporad2015quadratic}, and both having $M_k$ and $d_k$ in $A_k$, in contrast to only $M_k$ as in the proposed method, leads to some undesirable properties. 
 
First of, elements in $M_k$ and $d_k$ can be of different magnitude which can lead to numerical problems. This was partly resolved in \cite{bemporad2017numerically} by introducing a scaling factor $\beta$ at the cost of some extra overhead. 

Secondly, and more critically, $M_k$ depends on $H$ and $A$ while $d_k$ also depends on $f$ and $b$. This is of importance when numerical stability of the methods are improved with outer proximal-point iterations \cite{bemporad2017numerically} and when the QPs are solved in the context of linear MPC. In both of these applications, $H$ and $A$ remain constant while only $f$ and/or $b$ change between QP instances, which means that the computational complexity can be significantly reduced by starting the solver with the previous solution and reusing the $L$ and $D$ factors. 
For the proposed method this is straightforward since $L$ and $D$ are related to $M$ which, in turn, is related to the unchanged $H$ and $A$. For QPNNLS, however, $L$ and $D$ are related to $M$ \textit{and} $d$, where $d$ changes between iterations. Hence, to be able to reuse $L$ and $D$ between QP instances in QPNNLS, two rank-one updates are necessary. These rank-one updates introduce additional complexity compared with the proposed method.

Finally, the NNLS reformulation used in QPNNLS introduces an additional scaling factor which depends on the current iterate which is not necessary in the proposed method, leading to additional simplifications.

\subsection{Goldfarb-Idnani}
\label{ssec:gi}
Warm starting the proposed method is straightforward: given $\mathcal{W}_0$, any $\lambda_0$ satisfying $[\lambda_0]_{\mathcal{W}_0}\geq 0, [\lambda_0]_{\overline{\mathcal{W}}_0}=0$, suffices (however $[\lambda_0]_{\mathcal{W}_0}>0$ is preferable for numerical reasons). 
Warm starting is not as straightforward for the popular dual active-set method in \cite{10.1007BF02591962} (GI), which needs to be started in a nonnegative complementary basic solution, i.e., a nonnegative solution to \eqref{eq:linsystem}. Warm starting GI, hence, requires an additional procedure which finds a nonnegative complementary basic solution given an initial working set $\mathcal{W}_0$, increasing the computational burden. In fact, such a procedure is often similar to Steps \ref{step:innerstart}-\ref{step:innerstop} in Algorithm \ref{alg:dual-as}, i.e., is already embedded in the proposed method (cf, e.g., Alg. $5.3$ in \cite{bartlett2006qpschur}).

Another important difference between GI and the proposed algorithm is that GI, similar to many other active-set QP methods, computes a step direction from the current iterate $\lambda_k$ to $\lambda^*_k$ at each iteration, whereas Algorithm \ref{alg:dual-as} computes $\lambda^*_k$ directly. Some advantages of computing $\lambda^*_k$ directly are:

  (i) The availability of $\lambda^*_k$ allows for the update $\lambda_{k+1} \leftarrow \lambda^*_k$ in Algorithm \ref{alg:dual-as} when $\lambda^*_k\geq0$. This means that numerical errors in $\lambda_{k+1}$ only stem from solving \eqref{eq:linsystem} \textit{one} time. In contrast, by computing a step direction, numerical errors in the iterate accumulate each time a step is taken, i.e., there is no "resetting" mechanism for the numerical errors in the iterate in GI, as for Algorithm \ref{alg:dual-as}, when $\lambda^*_k \geq0$.  

  (ii) Less constraints can be classified as blocking when $\lambda^*_k$ is computed directly since only negative components of $\lambda^*_k$ can be blocking while, similarly, all negative components of the step direction are seen as possible blocking constraint for GI, where the former is always a subset of the latter. This ultimately leads to fewer computations for the former when removing constraints from $\mathcal{W}$, since fewer ratio tests (Step \ref{step:ratio} in Algorithm \ref{alg:dual-as}) have to be done. 

  (iii) In particular for GI, $\lambda^*_k$ can not be determined directly from the computed step direction since the step length along the step direction to reach $\lambda^*_k$ is unknown. Therefore, a primal iterate is also updated and used to determine which step length should be taken along the step direction from the iterate $\lambda_k$ to reach $\lambda^*_k$ (this is referred to the full step length in \cite{10.1007BF02591962}). 
	Performing these primal updates increases the computations needed compared with Algorithm \ref{alg:dual-as}, where only the dual iterate $\lambda_k$ needs to be updated in each iteration.

\subsection{Other active-set methods}
\label{ssec:other-as}
As is shown in \cite{best1984equivalence}, many active-set methods reported in the literature are mathematically equivalent, in the sense that they produce the same sequence of iterates before reaching the solution, and Algorithm \ref{alg:dual-as} belongs to the family of methods considered therein. For example, the proposed method is mathematically equivalent to the active-set algorithms presented in \cite{10.1007BF02591962,fletcher1971general,dantzig1998linear}. The differences between these active-set algorithms are numerical, e.g., how the systems of linear equations are solved and book-keeping of iterates.

The proposed method also shares the interpretation of the dual active-set methods in \cite{falt2019qpdas} (QPDAS) and in \cite{axehill2006mixed} (DRQP) as a primal active-set method applied to the dual problem \eqref{eq:qp-dual}. Instead of handling unbounded subproblems directly, as is done in the proposed method, QPDAS performs proximal-point iterations on the dual problem.  Furthermore, the factorization used for solving the subproblems is different.
DRQP differs from the proposed method in that it works on a sparse QP formulation rather than a dense one. By doing so, the subproblems are solved using a Riccati recursion, which leads to a linear computational complexity in the horizon of the MPC problem when solving the subproblems. DRQP cannot, however, detect primal infeasibility directly and does not perform low-rank updates to reduce computations, while the proposed algorithm does both (see Section~\ref{ssec:infeas} and Section~\ref{ssec:ldlupdate}, respectively, for details).

Recursively updating an LDL$^T$ factorization when solving the subproblems in an active-set method is described in the context of primal active-set methods in \cite{gill1974methods}. Since the constraints are particularly simple for the dual problem \eqref{eq:qp-dual}, some of the computations described in \cite{gill1974methods} simplify when the factorization is used in the context of the proposed dual active-set method. 

\section{Extensions}
\label{sec:extensions}
\subsection{Detecting infeasibility}
\label{ssec:infeas}
%A common way of solving MIQP is through branch and bound, where QPs are solved in a branch. When pruning such a branch, detecting infeasibility and obtain lower bounds for the QPs are of importance. Below we show how infeasibility and lower objective function bounds can be determined. 
  For a QP method to be reliable it needs to be able to detect if \eqref{eq:QP} has a primal feasible solution at all, i.e., if $\{x \in \mathbb{R}^n : A x\leq b\} \neq \emptyset$, otherwise the QP method might not be able to terminate in finite time for infeasible problems. Primal infeasibility can be detected in Step \ref{step:singblocking} in Algorithm \ref{alg:dual-as} if $\mathcal{B}=\emptyset$, as is shown by the following lemma.
\begin{lemma}
  \label{lem:unbound}
  Let $p_k$ satisfy $[p_k]_{\overline{\mathcal{W}}_k} = 0, M_k M_k^T [p_k]_{\mathcal{W}_k}=0$, and  $d^Tp_k<0$. Furthermore,  let $\mathcal{B}\triangleq\{i\in\mathcal{W}_k: [p_k]_i <0\}$. 
Then if $\mathcal{B}=\emptyset$, i.e., if $p_k \geq 0$, the QP in \eqref{eq:QP} is infeasible.
\end{lemma}
\begin{IEEEproof}
  Inserting $\alpha p_k$, for an arbitrary constant $\alpha >0$,  in the dual objective gives 
  \begin{equation*}
	\begin{aligned}
	  V(\alpha p_k)=& \alpha^2 \frac{1}{2}{p_k}^T M M^T p_k + \alpha d^T p_k  \\
	  =& \alpha^2 \frac{1}{2} [p_k]_{\mathcal{W}_k}^T M_k M_k^T [p_k]_{\mathcal{W}_k} + \alpha d^T p_k 
	  = \alpha d^T p_k,
	\end{aligned}
  \end{equation*}
  where the second equality follows from $[p_k]_{\overline{\mathcal{W}}_k} = 0$ and the last equality follows from $M_k M_k^T [p_k] _{\mathcal{W}_k}= 0$. Now, since $d^T p_k < 0$, $V(\alpha p_k)\to -\infty$ as $\alpha \to \infty$. Furthermore, $\alpha p_k$ is dual feasible since $\alpha p_k \geq 0, \:\forall \alpha >0$ ($p_k \geq 0$ by construction), making \eqref{eq:qp-dual} unbounded. 
The desired result follows from an unbounded dual problem being equivalent to an infeasible primal problem, see, e.g., \cite[Sec. 5.2.2]{boyd2004convex}.
\end{IEEEproof}
\subsection{Intermediary lower bounds on $J(x^*)$}
In some applications, for example when QP subproblems are solved as a part of solving mixed-integer quadratic programs (MIQPs) with branch-and-bound, having lower bound on $J(x^*)$ can reduce computations significantly \cite{fletcher1998numerical}. 

Since Algorithm \ref{alg:dual-as} operates on the dual problem, the well-known result from convex optimization that the dual function evaluated at a dual feasible point yields lower bounds on $J(x^*)$ (see, e.g., \cite[Sec. 5.1.3]{boyd2004convex}) can be used to efficiently compute such bounds. Concretely, we get the lower bound 
\begin{equation}
  \label{eq:low-bound}
  J(x^*) \geq \frac{1}{2}(\|M_k^T \lambda^*_k\|_2^2 - \|v\|_2^2) 
\end{equation}
 every time $\lambda^*_k \geq 0$ in Algorithm~\ref{alg:dual-as}.

 Moreover, by inserting $\lambda^*_k$ in \eqref{eq:qp-dual} and using \eqref{eq:linsystem} one gets that $J_d(\lambda^*_k) = -\tfrac{1}{2}\|M_k^T \lambda^*_k\|^2_2$, and, since Algorithm \ref{alg:dual-as} is a descent method, $\|M_k^T \lambda^*_k\|^2_2$ will increase in subsequent iterations, resulting in the lower bounds in \eqref{eq:low-bound} increasing (and becoming tight once $\mu_k \geq 0$).
\begin{remark}[Detecting cycling]
  $\|M_k^T \lambda^*_k\|^2_2$ can also be used to detect cycling in Algorithm \ref{alg:dual-as}, which can occur for ill-conditioned problems due to rounding errors. This can be done by checking whether $\|M_k^T \lambda^*_k\|$ increases every time $\lambda_k^* \geq 0$, which ensures that Algorithm \ref{alg:dual-as} is making progress in each iteration.
\end{remark}
\subsection{Bound constraints}
Often the constraints in QPs encountered in applications, for example in MPC, are given by upper and lower bounds in the form ${b^-\leq Ax\leq b^+}$. A naive way of handling these is to reformulate the constraints as $\tilde{A}\leq \tilde{b}$, with 
${\tilde{A} \triangleq[\begin{smallmatrix}
	A \\ -A 
\end{smallmatrix}]}, 
	\tilde{b} \triangleq[\begin{smallmatrix}
   b^+ \\-b^- 
	\end{smallmatrix}]$,
	which puts the QP in the form in \eqref{eq:QP}. The structure of the bound constraints can, however, be used to reduce the computational complexity and memory footprint of Algorithm \ref{alg:dual-as}, primarily when computing $\mu_k$. When exploiting the bound structure, each component of $[\lambda]_i$ corresponds to, instead of just the one-sided constraint $[A]_i x \leq [b]_i$, the two-sided constraint $[b^-]_i \leq [A]_ix \leq [b^+]_i$. To distinguish between whether the upper or lower bound is active, the sets $\mathcal{W}^+$ and $\mathcal{W}^-$, containing components corresponding to active upper and lower bounds, respectively, are introduced. Note that $\mathcal{W}^+ \cup \mathcal{W}^- = \mathcal{W}$.

To determine if optimality has been achieved or whether a constraint needs to be added to $\mathcal{W}$, we consider the primal slacks $\mu^{\text{+}}_k$ and $\mu^{\text{-}}_k$ for the upper and lower bounds, respectively, for the constraints not in $\mathcal{W}_k$, computed by
\begin{equation}
  \label{eq:mu-bounds}
  \begin{aligned}
	\:[\mu_k^{\text{+}}]_{\overline{\mathcal{W}}_k} &= \overline{M}_k M_k^T \lambda^*_k + [b^++Mv]_{\overline{\mathcal{W}}_k}, \\
	[\mu_k^{\text{-}}]_{\overline{\mathcal{W}}_k} &= -\overline{M}_k M_k^T \lambda^*_k - [b^-+Mv]_{\overline{\mathcal{W}}_k}. 
  \end{aligned}
\end{equation}

Hence, instead of $d$, the algorithm uses $d^+ \triangleq b^+ + Mv$ and $d^- \triangleq -(b^-+Mv)$.
Importantly, the relatively expensive matrix multiplication $\overline{M}_kM_k^T \lambda^*_k$ only has to be computed once in \eqref{eq:mu-bounds}, while it has to be computed twice if the naive formulation with $\tilde{A}$ and $\tilde{b}$ is used. Furthermore, if a component of $\mu_k^{\text{+}}$ is negative, the corresponding component of $\mu_k^{\text{-}}$ does not have to be computed since both the upper and lower bounds cannot be violated simultaneously (under the assumption that $b^- \leq b^+$, i.e. that the QP problem is not trivially infeasible). Optimality has been achieved if $\mu_k^{\text{+}}\geq 0$ and $\mu_k^{\text{-}}\geq 0$. Otherwise, the most negative component of $\mu_k^{\text{+}}$ or $\mu_k^{\text{-}}$ is added to $\mathcal{W}_k$.

When each component of $\lambda$ is the multiplier for both an upper and lower bound simultaneously, it does not have to be nonnegative anymore. Instead, components of $\lambda_k$ corresponding to active \textit{upper} bounds have to be \textit{nonnegative} while, conversely, components corresponding of active \textit{lower} bounds have to be \textit{nonpositive}. The condition $\lambda^*_k \geq 0$ is, hence, replaced by 
\begin{equation}
  [\lambda^*_k]_i \geq 0\:\:\forall i\in \mathcal{W}^+_k,\text{ and }[\lambda^*_k]_i \leq 0 \:\: \forall i\in \mathcal{W}^-_k. 
\end{equation}
Moreover, $\mathcal{B}$ is redefined as 
$\mathcal{B} \triangleq\{i\in \mathcal{W}^+:[\lambda^*_k]_i <0\} \cap\{i\in \mathcal{W}^{-}:[\lambda^*_k]_i >0\}$. For the singular case, $\mathcal{B}$ is redefined in a similar way but in terms of $p_k$ rather than $\lambda^*_k$.

Finally, when solving the subproblems $M_k M_k^T \lambda^*_k = -d_k$, the components of $d_k$ corresponding to active upper bounds should be replaced by elements of $d^+$ and, likewise, the components of $d_k$ corresponding to active lower bounds should be replaced by elements of $-d^-$.  
\begin{remark}[Box-constrained QP]
  When the constraints are in the simple form $b^{\text{-}}\leq x \leq b^{\text{+}}$, common in, e.g., MPC, then $M={R}^{-1}$, which is upper triangular. This additional structure can be exploited to reduce computations and the memory footprint further. 
\end{remark}
\subsection{Equality constraints}
\label{ssec:eqc}
Algorithm \ref{alg:dual-as} can also easily be extended to handle equality constraints in \eqref{eq:QP}. If the equality constraints are given as ${Gx=h}$, we can, similar to \eqref{eq:aux-def}, define 
\begin{equation}
  N \triangleq G R^{-1}, \quad w \triangleq R^{-1} f,\quad e \triangleq h + N w.
\end{equation}
The dual of the QP can then be stated as
\begin{equation}
  \begin{aligned}
	&\underset{\lambda\geq 0,\: \nu}{\text{minimize }} 
  \frac{1}{2}
  \begin{bmatrix}
   \nu \\
   \lambda
  \end{bmatrix}^T
  \begin{bmatrix}
   N \\ M  
  \end{bmatrix}
  \begin{bmatrix}
	N \\ M  
  \end{bmatrix}^T
  \begin{bmatrix}
   \nu\\
   \lambda
  \end{bmatrix}+
  \begin{bmatrix}
	e \\ d
  \end{bmatrix}^T
  \begin{bmatrix}
   \nu\\
   \lambda
  \end{bmatrix}
  \end{aligned}
\end{equation}
Essentially, equality constraints can be interpreted to always be active, i.e., be treated as being in $\mathcal{W}_k$ for all $k$. The corresponding modifications to Algorithm \ref{alg:dual-as} are, hence, to replace \eqref{eq:linsystem}, \eqref{eq:singdireq} and \eqref{eq:primal-slack} with 
\begin{subequations}
\begin{align}
  \begin{bmatrix}
   N \\M_k  
  \end{bmatrix}
  \begin{bmatrix}
   N \\ M_k 
  \end{bmatrix}^T 
  \begin{bmatrix}
	\nu^*_k \\ [\lambda^*_k]_{\mathcal{W}_k}
  \end{bmatrix}
  = 
  -\begin{bmatrix}
	e
	\\d_k
\end{bmatrix},\\
  \begin{bmatrix}
  N\\ M_k  
  \end{bmatrix}
  \begin{bmatrix}
   N \\ M_k 
  \end{bmatrix}^T 
  \begin{bmatrix}
	\tilde{\nu}_k \\ [p_k]_{\mathcal{W}_k}
  \end{bmatrix}
  = 0,\quad 
  d^T p_k
  <0,\\
  \mu_k = \overline{M}_k \begin{bmatrix}
   N\\M_k 
  \end{bmatrix}^T
  \begin{bmatrix}
	\nu^*_k \\ [\lambda^*_k]_{\mathcal{W}_k}
  \end{bmatrix}
  + \overline{d}_k, 
\end{align}
\end{subequations}
respectively.

Finally, factors $L$ and $D$ such that $N N^T = LDL^T$ are computed in the start of the algorithm. 
\subsection{Proximal-point iterations}
\label{ssec:proximal}
The numerics of any QP solver can be improved by performing outer proximal-point iterations, which results in a sequence of better conditioned QPs being solved. In particular, proximal-point iterations are given by
\begin{equation}
  \label{eq:QP-prox}
  \begin{aligned}
	x_{k+1} = &\underset{x}{\text{ argmin}}&&\frac{1}{2}x^T (H+\epsilon I) x + (f-\epsilon x_k)^T x \\
			  &\:\:\:\:\text{s.t.} &&Ax \leq b,\\
  \end{aligned}
\end{equation}
where $\epsilon>0$ is a regularization parameter. It can be shown (cf. \cite[Corollary 1]{bemporad2017numerically}) that $\lim_{k \to \infty} x_k = x^*$ when  \eqref{eq:QP-prox} is iteratively applied. 
As was shown in \cite{bemporad2017numerically}, combining outer proximal-point iterations with an active-set algorithm can lead to an efficient and numerically stable solver. The numerical stability of the proposed method can, hence, be improved by amending it with outer proximal-point iterations, summarized in Algorithm \ref{alg:proximal}.  In Algorithm \ref{alg:proximal}, $R_{\epsilon}$ is an upper Cholesky factor to $H+\epsilon I$, $M_{\epsilon} \triangleq A R_{\epsilon}$, $\eta>0$ is a tolerance for determining if a fixed point has approximately been reached, and \textbf{DAQP} refers to Algorithm \ref{alg:dual-as}. Also note that \textbf{DAQP} is heavily warm-started in Step \ref{step:daqp-sub} when each perturbed QP in the form \eqref{eq:QP-prox} is solved, reducing the computational burden significantly. 
\begin{algorithm}[H]
  \caption{Proximal-point iterations.}
  \label{alg:proximal}
  \begin{algorithmic}[1]
	\Require $\epsilon>0, M_{\epsilon},b,f, R_{\epsilon}^{-1}, \mathcal{W}, \lambda, x$
	\Ensure $x^*, \lambda^*, \mathcal{A}^*$
	\While{true}
	\State $v \leftarrow  R_{\epsilon}^{-T}(f-\epsilon x)$; \quad $d \leftarrow b + M_{\epsilon}v$ 
	\State $x_\text{old} \leftarrow x$
	\State $[x,\lambda,\mathcal{W}] \leftarrow \textbf{DAQP}(M_{\epsilon},d,v,R^{-1}_{\epsilon},\lambda,\mathcal{W})$\label{step:daqp-sub}
	\If{$\|x-x_{\text{old}}\|<\eta$}
	\State \textbf{return} $x, \lambda, \mathcal{W}$ 
	\EndIf
	\EndWhile
	%\State $k \leftarrow k+1$
  \end{algorithmic}
\end{algorithm}
%\begin{remark}
%  Since $L$ and $D$ used in Algorithm \ref{alg:dual-as} only depends on $M$, which does not change between the iterations, $L$ and $D$ from the previous iteration can \textit{directly} be used when Algorithm \ref{alg:dual-as} is warm started with $\mathcal{W}$ from the previous outer iteration. This is not the case for the method presented in \cite{bemporad2017numerically}, where $L$ and $D$ have to undergo low-rank updates between the outer iterations since its $L$ and $D$ factors also depends on $d$, which changes between the outer iterations.
%\end{remark}
\subsection{Exact complexity certification}
As was mentioned in Section \ref{sec:qpalg}, Algorithm \ref{alg:dual-as} can be interpreted as the active-set algorithm considered in \cite{journal2020}. Therein, a complexity certification method is proposed, which exactly determines the computational complexity for this active-set algorithm when QPs originating from a given multi-parametric quadratic program are to be solved. This complexity certification method can, hence, be used to determine the exact computational complexity of Algorithm \ref{alg:dual-as}. By doing so, worst-case bounds on the number of iterations and/or floating operations of Algorithm \ref{alg:dual-as} can be determined \textit{before} it is used in, e.g., an embedded MPC application.

\section{Numerical Experiments}
\label{sec:result}
\subsection{Comparison with QPNNLS}
\label{ssec:resqpnnls}
Numerical stability of the proposed method is compared with QPNNLS on a set of randomly generated small-scale QPs with varying condition numbers $\kappa(H)$. For each $\kappa(H)$, 100 QPs of size $n=25, m=100$ are generated and the worst-case distance from the optimal solution $x^*$ as well as the worst-case solution time are measured. 

We compare: The original formulation of QPNNLS presented in \cite{bemporad2015quadratic} (QPNNLS), the extended version of QPNNLS presented in \cite{bemporad2017numerically} in which numerical stability is improved by scaling and by performing outer proximal-point iterations (QPNNLS PROX), the proposed method given by Algorithm \ref{alg:dual-as} (DAQP), and, finally, Algorithm \ref{alg:dual-as} in conjuction with the outer proximal-point iterations given by Algorithm \ref{alg:proximal} (DAQP PROX). 

For all experiments, the primal feasibility tolerance is $\epsilon_p = 10^{-6}$ and, for the proximal-point iterations, $\epsilon=10^{-4}$ and $\eta=\sqrt{2^{-52}}\approx 1.5\cdot 10^{-8}$ (square root of machine epsilon for double precision). 
All of the methods are implemented in MATLAB using double precision and reference $x^*$ are obtained by using CPLEX with settings emphasizing numerical precision. Each QP is solved five times and the median execution time for these five runs is the reported solution time. 
\pgfplotstableread{data/wc_time2.dat}{\maxtime}
%\begin{tikzpicture}[scale=0.8]
%  \begin{axis}[xmode=log,xmin=100,xmax=10^10,ymin=2.05,ymax=6.5,
%	xlabel=$\kappa(H)$,
%	ylabel={Worst-case solution time [ms]},
%	legend style={at ={(0,1)},anchor=north west}, ymajorgrids,xmajorgrids,
%	x post scale=1.25,
%	legend cell align={left}]
%	\addplot [red,thick] table [x={cond}, y={daqp}] {\maxtime};
%	\addplot [black,thick] table [x={cond}, y={nnls}] {\maxtime};
%	\addplot [green,thick] table [x={cond}, y={daqp-prox}] {\maxtime};
%	\addplot [blue,thick] table [x={cond}, y={nnls-prox}] {\maxtime};
%	\legend{DAQP,QPNNLS,DAQP PROX,QPNNLS PROX};
%\end{axis}
%\end{tikzpicture}
\pgfplotstableread{data/wc_err2.dat}{\maxerror}
%\begin{tikzpicture}[scale=0.8]
%  \begin{axis}[ymode=log,xmode=log,xmin=100,xmax=10^10,ymin=10^-12,ymax=10^-4,
%	xlabel=$\kappa(H)$,ylabel={Worst-case $||x-x^*||_2$},ymajorgrids,xmajorgrids,
%	x post scale=1.25]
%	\addplot [black,very thick] table [x={cond}, y={nnls}] {\maxerror};
%	\addplot [red,very thick] table [x={cond}, y={daqp}] {\maxerror};
%	\addplot [blue,very thick] table [x={cond}, y={nnls-prox}] {\maxerror};
%	\addplot [green,very thick] table [x={cond}, y={daqp-prox}] {\maxerror};
%\end{axis}
%\end{tikzpicture}
\begin{figure*}[htpb] 
  \centering
  \subfloat[Worst-case solution time\label{subfig:soltime}]{%
	\begin{tikzpicture}[scale=0.8]
	  \begin{axis}[xmode=log,xmin=100,xmax=10^10,ymin=2.05,ymax=6.5,
		xlabel=$\kappa(H)$,
		ylabel={Worst-case solution time [ms]},
		legend style={at ={(0,1)},anchor=north west}, ymajorgrids,xmajorgrids,
		x post scale=1.25,
		legend cell align={left}]
		\addplot [set19c1,dashed,very thick] table [x={cond}, y={daqp}] {\maxtime};
		\addplot [set19c2,dashed,very thick] table [x={cond}, y={nnls}] {\maxtime};
		\addplot [set19c1,very thick] table [x={cond}, y={daqp-prox}] {\maxtime};
		\addplot [set19c2,very thick] table [x={cond}, y={nnls-prox}] {\maxtime};
		\legend{DAQP,QPNNLS,DAQP PROX,QPNNLS PROX};
	  \end{axis}
	\end{tikzpicture}
  }
  \hfill
  \subfloat[Worst-case distance to optimizer\label{subfig:dist}]{%
\begin{tikzpicture}[scale=0.8]
  \begin{axis}[ymode=log,xmode=log,xmin=100,xmax=10^10,ymin=10^-12,ymax=10^-4,
	xlabel=$\kappa(H)$,ylabel={Worst-case $||x-x^*||_2$},ymajorgrids,xmajorgrids,
	x post scale=1.25]
	\addplot [set19c1,dashed,ultra thick] table [x={cond}, y={daqp}] {\maxerror};
	\addplot [set19c2,dashed,ultra thick] table [x={cond}, y={nnls}] {\maxerror};
	\addplot [set19c2,ultra thick] table [x={cond}, y={nnls-prox}] {\maxerror};
	\addplot [set19c1,ultra thick] table [x={cond}, y={daqp-prox}] {\maxerror};
\end{axis}
\end{tikzpicture}
}
\caption{Comparison of DAQP, QPNNLS, with and without performing outer proximal-point iterations, on QPs with varying condition number $\kappa(H)$. The worst case solution time and distance to the optimizer reported for 100 randomly generated QPs with $n=25$ and $m=100$, for each $\kappa(H)$. 
 % Since solution times are based on MATLAB implementations, it is the relative, rather than the absolute, times that are of interest.
}
\label{fig:qpnnlsfig} 
\end{figure*}
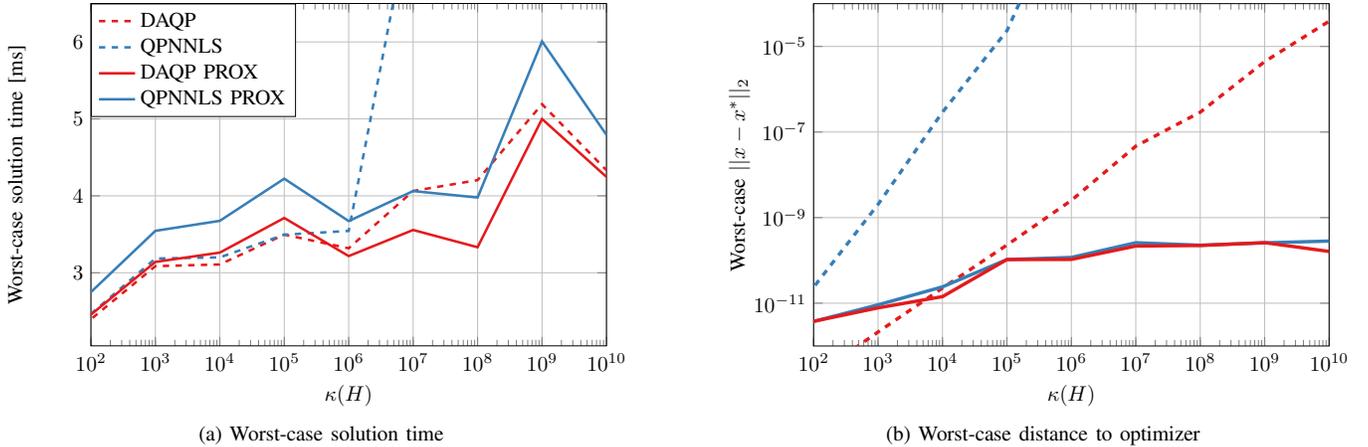
%\begin{tikzpicture}[scale=1]
%  \begin{axis}[ymode=log,minor tick num=1,xmin=10,xmax=50,
%	xlabel=$N$,ylabel={Worst-case solution time [s]},
%	legend style={at ={(0,1)},anchor=north west}]
%	\addplot [black,very thick] table [x={N}, y={prop}] {\maxtime};
%	\addplot [red,very thick] table [x={N}, y={qpnnls}] {\maxtime};
%	\addplot [blue,very thick] table [x={N}, y={goldfarb}] {\maxtime};
%	\legend{DAQP,QPNNLS,GI};
%\end{axis}
%\end{tikzpicture}
%\pgfplotstableread{data/max_error.dat}{\maxerror}
%\begin{tikzpicture}[scale=1]
%  \begin{axis}[ymode=log,minor tick num=1,xmin=10,xmax=50,
%	xlabel=$N$,ylabel={Worst-case $||x-x^*||_2$}]
%	\addplot [black,very thick] table [x={N}, y={prop}] {\maxerror};
%	\addplot [red,very thick] table [x={N}, y={qpnnls}] {\maxerror};
%	\addplot [blue,very thick] table [x={N}, y={goldfarb}] {\maxerror};
%\end{axis}
%\end{tikzpicture}

Worst-case results are shown Figure \ref{fig:qpnnlsfig}, with DAQP showing better numerical properties and worst-case solution time compared with QPNNLS, as is to be expected from the discussion in Section \ref{ssec:qpnnlscomp}.
Furthermore, Figure \ref{fig:qpnnlsfig} also illustrates that the robust numerical properties of QPNNLS PROX reported in \cite{bemporad2017numerically} extend to when DAQP is amended with proximal-point iterations. The worst-case solution time for DAQP PROX is, however, less than that of QPNNLS PROX, mainly because DAQP PROX can directly reuse the LDL$^T$ factorization between the outer proximal-point iterations, while QPNNLS PROX has to perform low-rank updates before reusing the factors.
\begin{remark}
  For some of the generated QPs with $\kappa(H)>10^6$, QPNNLS reached the iteration limit (which was set to 250) and was, hence, unable to return a solution within the primal infeasibility tolerance. This explains the jump in solution time for $\kappa(H)>10^6$ in Figure \ref{subfig:soltime}.
\end{remark}
\subsection{Model predictive control application}
The proposed algorithm is also tested for MPC of an ATFI-16 aircraft \cite{bemporad1997nonlinear} with varying prediction horizon $N$, which is a tutorial problem in the MATLAB Model Predictive Control Toolboox and was also considered in \cite{bemporad2015quadratic,bemporad2017numerically}. The system consists of two inputs, flaperon and elevator angles, which are upper and lower constrained, and two outputs, attack and pitch angles, with upper and lower limits imposed on the pitch angle. The system is simulated for 10 seconds with an angle of attack reference of $0$ and a pitch angle reference shifting between $\pm10^\circ$. For prediction horizon $N$, the resulting QPs have dimensions $n= 2N+1$, $m = 4N+2(N-1)$, where the extra optimization variable is due to the output constraints being softened to always guarantee that a primal feasible solution exists.
Because of the system having unstable poles, the Hessian has a fairly high condition number ($\kappa(H) > 10^{10}$).  For the weights in the MPC formulation, the same settings as in the MATLAB Model Predictive Control Toolbox tutorial were used. 

Since the same MPC problem also was considered in \cite{bemporad2015quadratic,bemporad2017numerically} and that the proposed method is a direct improvement of the methods therein (motivated in Section \ref{ssec:qpnnlscomp} and highlighted by the experiments in Section \ref{ssec:resqpnnls}), the favourable results for QPNNLS reported therein immediately also holds for the proposed method, and are in fact even strengthened. Nevertheless, we compare the worst-case solution time for the proposed method with some additional QP solvers which are used in the context of MPC. Herein, we focus on the case when the MPC problem is reformulated as a dense QP. 
\begin{remark}[Exploiting sparisty]
  We want to stress that our focus here is on solving dense QPs from the condensed MPC problem formulation. When solving large problems, other solvers that use sparse formulations \cite{frison2020hpipm,stellato2020osqp,axehill2006mixed,quirynen2020presas,frasch2015parallel} might be more efficient. If the ideas herein can be modified to exploit sparsity is a topic for future research. 
\end{remark}

We compare a C implementation of Algorithm \ref{alg:dual-as} (DAQP) and Algorithm \ref{alg:proximal} (DAQP PROX) with: (i) The parametric active-set method presented in \cite{ferreau2014qpoases} (qpOASES\_e), coded in C; (ii) The operator splitting method presented in \cite{stellato2020osqp} (OSQP), coded in C; (iii) The interior-point method presented in \cite{frison2020hpipm} (HPIPM), coded in C; (iv) The dual active-set method presented in \cite{schmid1994quadratic} (QPKWIK) which is based on \cite{10.1007BF02591962}, with generated C code from MATLAB. 

To get reliable solution times, each QP was solved $15$ times and the median solution time was used as the reported solution time for the QP. When possible and relevant, the solvers were provided with precomputed matrix factorizations, e.g., the Cholesky factor of $H$ was precomputed for DAQP, qpOASES\_e and QPKWIK and only solve time was considered for OSQP. 

\pgfplotstableread{data/MPC/daqp_MPC.dat}{\daqpMPC}
\pgfplotstableread{data/MPC/daqp_prox_MPC.dat}{\daqpproxMPC}
\pgfplotstableread{data/MPC/qpOASES_MPC.dat}{\qpOASESMPC}
\pgfplotstableread{data/MPC/QPKWIK_MPC.dat}{\qpkwikMPC}
\pgfplotstableread{data/MPC/OSQP_MPC.dat}{\osqpMPC}
\pgfplotstableread{data/MPC/HPIPM.dat}{\hpipmMPC}

\pgfplotstableread{data/MPC/Jdiffs.dat}{\Jdiffs}
%\begin{figure} 
%	\begin{tikzpicture}[scale=0.8]
%	  \begin{axis}[xmin=5,xmax=30,
%		ymin=-5,ymax=6,
%		xlabel={Prediction horizon ($N$)},
%		ytick={-4,-2,0,2,4},
%		yticklabels={$-10^{-6}$,-$10^{-7}$,$0$,$10^{-7}$,$10^{-6}$},
%		ylabel={$\texttt{mean}\big[J(x^*)-J(x^*_{\text{DAQP}})\big]$},
%		legend style={at ={(1,0)},anchor=south east}, ymajorgrids,yminorgrids,xmajorgrids,
%		x post scale=1.25,
%		legend cell align={left},legend columns=2,
%		]
%		\addplot [set19c1,very thick, dashed] table [x={horizon}, y={daqp}] {\Jdiffs};
%		\addplot [set19c2,very thick] table [x={horizon}, y={daqpprox}] {\Jdiffs};
%		\addplot [set19c3,very thick] table [x={horizon}, y={osqp}] {\Jdiffs};
%		\addplot [set19c4,very thick] table [x={horizon}, y={hpipm}] {\Jdiffs};
%		\addplot [set19c5,very thick] table [x={horizon}, y={qpkwik}] {\Jdiffs};
%		\addplot [yellow!75!black,very thick] table [x={horizon}, y={qpoases}] {\Jdiffs};
%		%\legend{DAQP PROX,OSQP, HPIPM, QPKWIK, qpOASES\_e};
%	  \end{axis}
%	\end{tikzpicture}
%	\caption{Worst-case solution time when solving the QPs arising during the MPC of an ATFI-16 aircraft in simulation, for varying prediction horizons $N$. The resulting QPs have dimensions $n=2N+1$, $m=4N+2(N-1)$. All of the solvers were exectued on an Intel 2.7 GHz i7-7500U CPU.}
%\label{fig:mpcfig} 
%\end{figure}
\begin{figure*}[htpb] 
  \centering
  \subfloat[Worst-case solution time\label{subfig:mpcsoltime}]{%
	\begin{tikzpicture}[scale=0.8]
	  \begin{axis}[xmin=5,xmax=30,
		ymode=log,ymin=10^-3,ymax=2*10^1,
		xlabel={Prediction horizon ($N$)},
		ylabel={Worst-case solution time [ms]},
		legend style={at ={(1,0)},anchor=south east}, ymajorgrids,yminorgrids,xmajorgrids,
		x post scale=1.25,
		legend cell align={left},legend columns=2,
		]
		\addplot [set19c1,very thick] table [x={horizon}, y={wc}] {\daqpMPC};
		\addplot [set19c2,very thick] table [x={horizon}, y={wc}] {\daqpproxMPC};
		\addplot [set19c3,very thick] table [x={horizon}, y={wc}] {\osqpMPC};
		\addplot [set19c4,very thick] table [x={horizon}, y={wc}] {\hpipmMPC};
		\addplot [set19c5,very thick] table [x={horizon}, y={wc}] {\qpkwikMPC};
		\addplot [yellow!75!black,very thick] table [x={horizon}, y={wc}] {\qpOASESMPC};
		\legend{DAQP, DAQP PROX,OSQP, HPIPM, QPKWIK, qpOASES\_e};
	  \end{axis}
	\end{tikzpicture}
  }
  \hfill
  \subfloat[Average cost difference to the solution of DAQP\label{subfig:mpc-quality}]{%
	\begin{tikzpicture}[scale=0.8]
	  \begin{axis}[xmin=5,xmax=30,
		ymin=-5,ymax=6,
		xlabel={Prediction horizon ($N$)},
		ytick={-4,-2,0,2,4},
		yticklabels={$-10^{-5}$,$-10^{-7}$,$0$,$10^{-7}$,$10^{-5}$},
		minor y tick num=1,
		ylabel={$\texttt{mean}\big[J(x^*)-J(x^*_{\text{DAQP}})\big]$},
		legend style={at ={(1,0)},anchor=south east}, ymajorgrids,yminorgrids,xmajorgrids,
		x post scale=1.25,
		legend cell align={left},legend columns=2,
		]
		\addplot [set19c1,very thick, dashed] table [x={horizon}, y={daqp}] {\Jdiffs};
		\addplot [set19c2,very thick] table [x={horizon}, y={daqpprox}] {\Jdiffs};
		\addplot [set19c3,very thick] table [x={horizon}, y={osqp}] {\Jdiffs};
		\addplot [set19c4,very thick] table [x={horizon}, y={hpipm}] {\Jdiffs};
		\addplot [set19c5,very thick] table [x={horizon}, y={qpkwik}] {\Jdiffs};
		\addplot [yellow!75!black,very thick] table [x={horizon}, y={qpoases}] {\Jdiffs};
		%\legend{DAQP PROX,OSQP, HPIPM, QPKWIK, qpOASES\_e};
	  \end{axis}
	\end{tikzpicture}
}
\caption{Worst-case solution time and average solution quality when solving QPs encountered during the MPC of an ATFI-16 aircraft in simulation for varying prediction horizons $N$. The resulting QPs have dimensions $n=2N+1$, $m=4N+2(N-1)$. All of the solvers were exectued on an Intel 2.7 GHz i7-7500U CPU.}
\label{fig:mpcfig} 
\end{figure*}
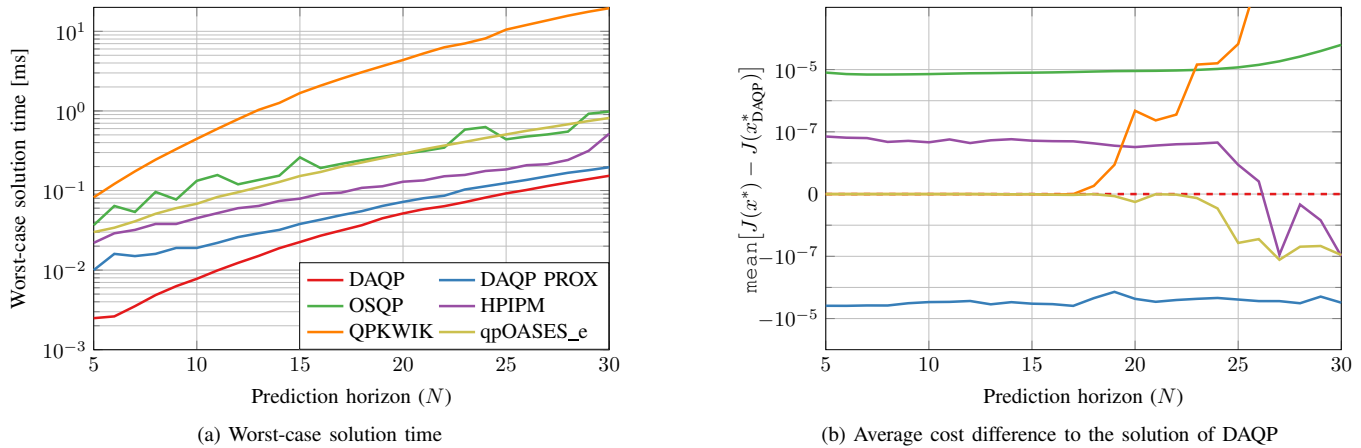

The worst-case solution times are reported in Figure \ref{subfig:mpcsoltime}, where DAQP and DAQP PROX outperform the other solvers. Note, however, that HPIPM and OSQP scale better with $N$ and will, as $N$ grows larger, sooner or later outperform active-set methods. Still, DAQP seems superior on small to medium size MPC problems, which is a common scope of problems where embedded model predictive control is employed. 
Moreover, the average quality of the solutions computed by the solvers are compared with the solution of DAQP in Figure \ref{subfig:mpc-quality} by considering differences in the objective function. A positive value, hence, means that the solution is on average worse than that of DAQP. Importantly, Figure \ref{subfig:mpc-quality} illustrates that in terms of solution quality, which is correlated with numerical robustness, DAQP PROX outperforms the other solvers. 
\section{Conclusion}
In this paper we have presented a dual active-set solver which is efficient, reliable, and simple to code, making it suitable for use in real-time model predictive control applications. Numerical experiments show that the proposed method can outperform state-of-the-art QP algorithms for QPs encountered in embedded MPC applications, both in terms of computational complexity and numerical robustness. 

Future work includes combining the ideas herein with low-rank updates of the Riccati factorization \cite{nielsen2017low}, to investigate if these ideas can lead to an efficient solver for problems of larger size. 
\bibliographystyle{IEEEtran}
\linespread{0.95}\selectfont
\bibliography{lib.bib}
\end{document}